\documentclass[12pt,a4paper]{amsart}
\usepackage[colorlinks,pagebackref,hypertexnames=false]{hyperref}
\usepackage[a4paper,left=1.8cm,right=1.8cm,top=2.5cm,bottom=2.5cm]{geometry}
\usepackage{lmodern}
\usepackage{bbm}
\usepackage{float}
\usepackage{xurl}
\usepackage{ragged2e}
\usepackage{enumitem}
\setlist[itemize]{noitemsep, topsep=0pt, partopsep=0pt}
\setlist[enumerate]{noitemsep, topsep=0pt, partopsep=0pt}
\setlist[description]{noitemsep, topsep=0pt, partopsep=0pt}

\usepackage{etoolbox}
\AtBeginEnvironment{thebibliography}{%
  \setlength{\itemsep}{0pt}%
  \setlength{\parskip}{0pt}%
  \setlength{\topsep}{0pt}%
  \setlength{\partopsep}{0pt}%
}

\newlist{defenum}{enumerate}{3}
\setlist[defenum]{label=(\alph*),leftmargin=*,align=left}
\newlist{defsubenum}{enumerate}{1}
\setlist[defsubenum]{label=(\roman*),leftmargin=*,align=left}

\hypersetup{pdftitle={An atomic approach to Shokurov's conic-bundle conjecture}}
\newtheorem{theorem}{Theorem}[section]
\newtheorem{proposition}[theorem]{Proposition}
\newtheorem{lemma}[theorem]{Lemma}
\newtheorem{corollary}[theorem]{Corollary}
\newtheorem{maintheorem}{Theorem}

\newtheorem{maincorollary}[maintheorem]{Corollary}
\theoremstyle{remark}
\newtheorem{remark}[theorem]{Remark}

\title{An atomic approach to Shokurov's conic-bundle conjecture}

\begin{document}

\author[Cavenaghi]{Leonardo F. Cavenaghi}
\address{Institute of Mathematics and Informatics, Bulgarian Academy of Sciences}
\email{leonardofcavenaghi@gmail.com}

\author[Katzarkov]{Ludmil Katzarkov}
\address{Simons Center for Geometry and Physics, Stony Brook University, Stony Brook, NY 11794, USA. \& International Center for Mathematical Sciences (ICMS), Sofia, Bulgaria. \& Institute of Mathematics and Informatics, Bulgarian Academy of Sciences, Sofia, Bulgaria.} 
\email{lkatzarkov@gmail.com}

\author[Pomerleano]{Daniel Pomerleano}
\address{UMASS BOSTON, 100 Morrissey Blvd., Boston, MA 02125} 
\email{daniel.pomerleano@umb.edu}

\author[Zhang]{Zhijia Zhang}
\address{Department of Mathematics, Columbia University, New York, NY 10027,
USA} 
\email{zhang.zhijia@columbia.edu}

\begin{abstract}
In this paper, we apply a weaker version of Hodge atoms to show that a standard conic bundle over a smooth complex projective rational surface $S$ is irrational if $|2K_S+\Delta|\ne\emptyset$, where $\Delta$ is the discriminant divisor. This proves Shokurov's conjecture on rationality of conic bundle threefolds. As an application, we deduce precisely which nodal quartic double solids are rational. 
\end{abstract}
\maketitle

\section{Introduction}
\label{sec:introduction}

Conic bundles occur naturally in the birational classification of complex threefolds.\footnote{In this paper, all varieties are understood to be over $\mathbb{C}$.} Their rationality problem relates the geometry of the discriminant curve to the intermediate
Jacobian of the total space; see \cite{Prokhorov} for an account of this theory. Inspired by the recently introduced theory of Hodge atoms of Katzarkov--Kontsevich--Pantev--Yu \cite{KKPY}, we study this problem through the algebraic structure of big quantum cohomology.

A \emph{standard conic bundle} is a surjective projective flat morphism $\pi: X\to S$ from a smooth projective threefold to a smooth projective surface, with plane conics as fibers and relative Picard number $\rho(X/S)=1$. Its discriminant $\Delta\subset S$ is the reduced curve over which the fibers are singular; when nonempty, it has at worst ordinary double points. We assume throughout that $S$ is rational and write
\begin{equation}\label{eq:adjoint}
\Delta_S:=2K_S+\Delta.
\end{equation}

Shokurov conjectured that the effectivity of $\Delta_S$ obstructs rationality \cite[Conjecture~10.3]{Shokurov}; see also \cite[Conjecture~2.2]{CPS2019}, and an equivalent conjecture in \cite{Isk1}. Here we prove this implication:

\begin{maintheorem}\label{thm:conic-bundle-effective-intro}
Let $\pi: X\to S$ be a projective standard conic bundle over a smooth projective rational surface, with discriminant $\Delta$. If
\[
|\Delta_S|\neq\emptyset,
\]
then $X$ is irrational.
\end{maintheorem}

Theorem~\ref{thm:conic-bundle-effective-intro} is first proven in the case where $\Delta_S$ is nef. Then, we use birational modifications of standard conic bundles to reduce the effective case to this one. Next, let us explain how a quantum-cohomology computation plays the main role in our argument.

\subsection{The quantum calculation}
 
In \cite{KKPY}, Katzarkov--Kontsevich--Pantev--Yu combine Hodge theory with Gromov--Witten theory to construct Hodge atoms. Their construction relates the spectral decomposition of quantum
multiplication by the Euler vector field, along the locus of Hodge classes, to the Hodge-theoretic structures carried by its summands. The behavior of these objects under blowups provides obstructions to rationality.

For our purposes, we use the formal algebraic approach of \cite{cai2026quarticthreefoldsymplecticallyirrational}. We study
quantum multiplication $\boldsymbol\kappa_{t'}=\mathsf{Eu}_{t'}\star_{t'}(-)$ on the even big quantum cohomology algebra over a specific field $\mathbb K_X$ defined in Section~\ref{sec:quantum}. The relevant obstruction is that, for a smooth projective rational threefold, every eigenvalue of this operator has algebraic multiplicity at most two \cite[Theorem~5.6]{cai2026quarticthreefoldsymplecticallyirrational}. Here we produce an eigenvalue of larger multiplicity by constructing a local Artinian factor:

\begin{maintheorem}\label{thm:conic-bundle-local-factor-intro}
Let $\pi: X\to S$ be a projective standard conic bundle over a smooth projective rational surface. Assume
that $\Delta_S$ is nef, and put $\rho_S:=\operatorname{rank}\mathsf{NS}(S)$.
Then
\[
QH^{\mathrm{even}}(X)
\otimes_{\mathsf{Nov}^{\mathrm{red}}_X}\mathbb K_X
\]
has a local Artinian factor $\mathcal H_-$ satisfying
\[
\begin{aligned}
\dim_{\mathbb K_X}\mathcal H_-&=\rho_S+2,&
\dim_{\mathbb K_X}\operatorname{rad}(\mathcal H_-)&=\rho_S+1,\\
\dim_{\mathbb K_X}\bigl(\operatorname{rad}(\mathcal H_-)^2\bigr)&=1,&
\operatorname{rad}(\mathcal H_-)^3&=0.
\end{aligned}
\]
The restriction of $\boldsymbol\kappa_{t'}$ to $\mathcal H_-$ has a single eigenvalue, whose algebraic multiplicity on the full even quantum cohomology is exactly $\rho_S+2$. In particular, $X$ is
irrational.
\end{maintheorem}

The proof of Theorem~\ref{thm:conic-bundle-local-factor-intro} has two steps. First, we compute the \emph{vertical quantum product}, obtained by retaining only curve classes contracted by $\pi$. This is loosely similar to the ``vertical quantum cohomology" computation that appears in \cite[\S 5]{Iritani}. After extending coefficients and inverting $q^\ell$, where $\ell$ is the class of a line component of a singular fiber, this algebra splits into two copies of $H^{\mathrm{even}}(S,\mathbb Q)$ with its cup product.

Second, we include the contributions of nonvertical curves and lift this decomposition to the full quantum product. The main point is to prove that one of the two lifted summands remains local over the algebraically closed fraction field. Nefness of $\Delta_S$ constrains the nonvertical contributions and gives the radical powers stated in the theorem.
Finally, we show that the eigenvalue of $\boldsymbol\kappa_{t'}$ on this local factor does not occur on the other summand, so its algebraic multiplicity is exactly $\rho_S+2$.

\subsection{Birational consequences}

Theorem~\ref{thm:conic-bundle-effective-intro}, together with the known results, gives the rationality criterion formulated in \cite[Conjecture~1.2]{Prokhorov}. Here $J(X)$ denotes the intermediate Jacobian with its canonical principal polarization. We prove:

\begin{maincorollary}\label{cor:conic-rationality-criterion-intro}
Let $\pi\colon X\to S$ be a standard conic bundle over a smooth projective rational surface, with discriminant $\Delta$. Then $X$ is rational if and only if $|2K_S+\Delta|=\emptyset$ and the intermediate Jacobian $J(X)$ is isomorphic, as a principally polarized abelian variety, to a product of Jacobians of smooth projective curves.
\end{maincorollary}

As a last contribution, we apply the calculation to \emph{nodal quartic double solids}, that is, double covers of $\mathbb P^3$ branched over quartic surfaces whose only singularities are (at worst) ordinary double points. Recall that a normal variety is \emph{$\mathbb Q$-factorial} if every Weil divisor has a positive multiple that is Cartier.

\begin{maincorollary}\label{thm:quartic-rationality-classification-intro}
Let $Y$ be a nodal quartic double solid with $N$ nodes. Then
\[
Y\text{ is rational}
\quad\Longleftrightarrow\quad
Y\text{ is not }\mathbb Q\text{-factorial and }N\neq6.
\]
\end{maincorollary}

Cheltsov--Przyjalkowski--Shramov established the connection between nodal quartic double solids and Shokurov's conjecture in \cite[Corollary~4.28]{CPS2019}.
For a $\mathbb Q$-factorial double solid with $N\geq1$ nodes, its standard model has $\Delta\sim-2K_S$ and $\rho_S=N$. Theorem~\ref{thm:conic-bundle-local-factor-intro} therefore gives a local factor of dimension $N+2$ on this smooth model. The non-$\mathbb Q$-factorial case and the six-node exception are given by \cite[Theorems~1.2 and~1.5]{CPS2019}. Moreover, Clemens' criterion expresses $\mathbb Q$-factoriality as independence of the conditions imposed by the nodes on quadrics in $\mathbb P^3$ \cite[Theorem~1.6]{CPS2019}. Building on that, in Section~\ref{sec:nodal-quartic-double-solids}, we give a rationality test in terms of a matrix-rank condition.

These conclusions fit into a long history of irrationality results for quartic double solids. Artin--Mumford constructed a ten-nodal example
whose smooth resolution has nonzero torsion in $H^3(-,\mathbb Z)$ and is therefore not stably rational
\cite{ArtinMumford}. Beauville's Prym-variety approach gives the one-node irrationality result and explains the limitation of the intermediate-Jacobian obstruction for examples with at least seven nodes \cite[Example~4.10.4]{BeauvillePrym}. For smooth quartic double solids, the intermediate-Jacobian approach was further developed by Tikhomirov \cite{Tikhomirov81} and Voisin \cite{Voisin88}. More recently, using decomposition of the diagonal, Voisin proved stable irrationality for very general double solids with at most seven nodes in the families specified
in \cite[Theorem~1.1 and Remark~1.2]{Voisin}. Achter--Casalaina-Martin--Vial extended this to cases with at most nine nodes \cite[Theorem~3 and Section~17]{ACV}. 

This paper is the starting point for a new direction, patching together the minimal model program and the theory of Hodge atoms; further developments shall appear in \cite{CKK-lectures, BMKG}, among others.

\subsection{This paper is organized as follows}
In Section~\ref{sec:quantum} we recall the basics of quantum cohomology and fix our conventions. Further technicalities used in subsequent sections are also discussed. In particular, the irrationality criteria used here are recalled. In Section~\ref{sec:conic-bundle-quantum} we compute the vertical product, extend its decomposition, and prove Theorem~\ref{thm:conic-bundle-local-factor-intro}. In Section~\ref{sec:conic-bundle-birational-reduction} we prove Theorem~\ref{thm:conic-bundle-effective-intro} and Corollary~\ref{cor:conic-rationality-criterion-intro}. Finally, in Section~\ref{sec:nodal-quartic-double-solids} we prove Corollary~\ref{thm:quartic-rationality-classification-intro} and give the explicit criterion on the nodes.

\subsection*{AI methodology}
The argument for Lemma~\ref{lem:nef-nonsplitting} was developed through a series of interactions with ChatGPT 6. This argument (along with the rest of the text) was written up by the authors.   

\subsection*{Acknowledgements}

   L.~F.~Cavenaghi and L.~Katzarkov are supported by the Simons Foundation, grant SFI-MPS-T-Institutes-00007697, and the Ministry of Education and Science of the Republic of Bulgaria, grant DO1-239/10.12.2024. This work was done while the first-named author was visiting the Simons Center for Geometry and Physics, for which he thanks the hospitality and nice work environment. L.~Katzarkov is supported by the Simons Investigators Award (no.003136), NSF FRG grant DMS-2245099, and NSF FRG grant DMS-2255039. D. Pomerleano was partially supported by NSF grant DMS-2604648. The authors thank Yuri Tschinkel for constant support and interest.

\section{Basics of Quantum Cohomology}
\label{sec:quantum}
Here, we recall the basic features of quantum cohomology needed to prove our results. Let $X$ be a smooth complex projective variety. We denote by
\[
N_1(X,\mathbb Z)
:=
\operatorname{im}\left(
\mathsf{CH}_1(X)
\xrightarrow{\operatorname{cl}}
H_2(X,\mathbb Q)
\right)
\]
the subgroup of $H_2(X,\mathbb Q)$ consisting of classes of algebraic $1$-cycles with integral coefficients, and let $\mathsf{NE}(X,\mathbb Z)\subset N_1(X,\mathbb Z)$ be the monoid of effective curve classes, including the zero class. Since $X$ is projective, for any integral ample divisor $H_X$ on $X$ there are only finitely many effective curve classes $\beta\in\mathsf{NE}(X,\mathbb Z)$ of bounded degree $\int_\beta H_X\leq k$. In particular, each class $\beta\in\mathsf{NE}(X,\mathbb Z)$ admits only finitely many decompositions $\beta=\beta_1+\beta_2$ with $\beta_1,\beta_2\in\mathsf{NE}(X,\mathbb Z)$ \cite[Lemmas~19--20]{LeePandharipande}. We may therefore define the \emph{Novikov ring} $\mathcal N_X$ as the completed monoid algebra of formal series
\[
\mathcal N_X:=
\mathbb Q[[\mathsf{NE}(X,\mathbb Z)]]
=
\left\{
\sum_{\beta\in\mathsf{NE}(X,\mathbb Z)} a_\beta q^\beta
\ :\ a_\beta\in\mathbb Q
\right\},
\qquad
q^\beta q^{\beta'}=q^{\beta+\beta'},
\]
endowed with the complete separated linear topology defined by the degree filtration
\[
I_k:=
\left\{
\sum_{\beta}a_\beta q^\beta:
a_\beta=0\ \text{whenever}\ \int_\beta H_X<k
\right\},
\qquad k\geq 1.
\]
Since $N_1(X,\mathbb Z)$ is torsion-free, its group algebra is an integral domain. Taking the nonzero parts of least ample degree shows that $\mathcal N_X$ is an integral domain as well.
The formal power-series extensions used below are therefore integral domains as well.

Here and throughout, let $H^\ast(X,\mathbb Q)$ denote the rational Betti cohomology of $X$, and set
\[
H^{\mathrm{even}}(X,\mathbb Q)
:=\bigoplus_{k\geq 0}H^{2k}(X,\mathbb Q).
\]
Write $\rho:=\operatorname{rank}\mathsf{NS}(X)$ for the Picard number of $X$. Fix a homogeneous $\mathbb Q$-basis
\[
\{T_0=\mathbf 1,T_1,T_2,\ldots,T_\rho,
T_{\rho+1},\ldots,T_m\}
\]
of $H^{\mathrm{even}}(X,\mathbb Q)$ such that $T_1,\ldots, T_\rho\in H^2(X,\mathbb Q)$ are the rational cohomology classes corresponding, under the first Chern class map, to a $\mathbb Z$-basis of $\mathsf{NS}(X)/\operatorname{tors}$. By the Lefschetz $(1,1)$-theorem, their $\mathbb Q$-span is
\[
\operatorname{im}\left(
c_1:\mathsf{NS}(X)\otimes_{\mathbb Z}\mathbb Q
\longrightarrow H^2(X,\mathbb Q)
\right)
=
H^2(X,\mathbb Q)\cap H^{1,1}(X),
\]
where the intersection is taken inside $H^2(X,\mathbb C)$ via the natural map $H^2(X,\mathbb Q)\to H^2(X,\mathbb C)$. Intersection with the chosen divisor classes defines an injective
homomorphism
\[
N_1(X,\mathbb Z)\longrightarrow\mathbb Z^\rho,
\qquad
\beta\longmapsto
\left(\int_\beta T_1,\ldots,\int_\beta T_\rho\right),
\]
so that we use the Laurent-monomial notation
\[
q^\beta=\prod_{i=1}^{\rho}q_i^{\int_\beta T_i}.
\]
One remarks, however, that the symbols $q_i$ are auxiliary coordinates and need not themselves belong to $\mathcal N_X$.

Let $t_0,\ldots,t_m$ be the dual coordinates to $T_0,\ldots,T_m$, viewed as formal variables. We define the \emph{reduced big Novikov ring} by
\[
\mathsf{Nov}^{\mathrm{red}}_X:=\mathcal N_X[[(t_i)_{i\not\in\{1,\ldots,\rho\}}]],
\]
and set
\[
t':=\sum_{i\not\in\{1,\ldots,\rho\}} t_i T_i \in H^{\mathrm{even}}(X,\mathbb Q)\otimes_{\mathbb Q}\mathsf{Nov}_X^{\mathrm{red}}.
\]
The (even)\footnote{hereafter, we abandon the adjective ``even''} \emph{big quantum cohomology ring}, denoted by $QH^{\mathrm{even}}(X)$, is the free $\mathsf{Nov}_X^{\mathrm{red}}$-module
\[
H^{\mathrm{even}}(X,\mathbb Q)\otimes_{\mathbb Q}\mathsf{Nov}_X^{\mathrm{red}}
\]
equipped with the \emph{big quantum product} $\star_{t'}$, which is the $\mathsf{Nov}_X^{\mathrm{red}}$-bilinear product uniquely determined for $\gamma_1,\gamma_2,\gamma_3\in H^{\mathrm{even}}(X,\mathbb Q)$ by
\[
\left(\gamma_1\star_{t'}\gamma_2,\gamma_3\right)_X
:=
\sum_{\beta\in\mathsf{NE}(X,\mathbb Z)}\sum_{n\geq 0}
\frac{1}{n!}\prod_{i=1}^{\rho}q_i^{\int_{\beta}T_i}
\langle \gamma_1,\gamma_2,\gamma_3,\underbrace{t',\ldots,t'}_{n\text{ times}}\rangle_{0,n+3,\beta}.
\]
Here $(\alpha,\gamma)_X:=\int_X\alpha\smile\gamma$ is the Poincar\'e pairing, extended bilinearly over $\mathsf{Nov}^{\mathrm{red}}_X$. We use primary genus-zero Gromov--Witten invariants, extended multilinearly in their insertions. Writing $\pi_{\mathbb Q}:H_2(X,\mathbb Z)\to H_2(X,\mathbb Q)$, our convention is
\[
\langle\alpha_1,\ldots,\alpha_N\rangle_{0,N,\beta}
:=
\sum_{\substack{\widetilde\beta\in H_2(X,\mathbb Z)\\
                 \pi_{\mathbb Q}(\widetilde\beta)=\beta}}
\int_{[\overline{\mathcal M}_{0,N}(X,\widetilde\beta)]^{\mathrm{vir}}}
\prod_{j=1}^{N}\operatorname{ev}_j^*(\alpha_j),
\]
for $N\geq0$ when $\beta\neq0$, and for $N\geq3$ when $\beta=0$. The sum is finite because $\operatorname{Tor}H_2(X,\mathbb Z)$ is finite; empty moduli spaces contribute zero. We refer to
\cite{fulton1997notesstablemapsquantum,Kontsevich1994} for the basic axioms.

\subsection{Quantum multiplication as a multiplication on the logarithmic tangent bundle}

For $1\leq i\leq\rho$, let $q_i\partial_{q_i}$ be the continuous derivation defined by
\[
q_i\partial_{q_i}(q^\beta)=\left(\int_\beta T_i\right)q^\beta,
\qquad q_i\partial_{q_i}(t_j)=0.
\]
Consider the free module of logarithmic vector fields
\[
\mathcal T_X^{\log}
:=
\bigoplus_{i=1}^{\rho}\mathsf{Nov}^{\mathrm{red}}_X q_i\partial_{q_i}
\ \oplus\
\bigoplus_{i\notin\{1,\ldots,\rho\}}
\mathsf{Nov}^{\mathrm{red}}_X\partial_{t_i}.
\]
We identify it with $H^{\mathrm{even}}(X,\mathbb Q)\otimes_{\mathbb Q} \mathsf{Nov}^{\mathrm{red}}_X$ via
\[
q_i\partial_{q_i}\longleftrightarrow T_i \quad (1\leq i\leq\rho),
\qquad
\partial_{t_i}\longleftrightarrow T_i
\quad (i\notin\{1,\ldots,\rho\}).
\]
Write $d_i:=\deg T_i$ for the real cohomological degree of $T_i$, and assign the formal weights
\[
\deg q^\beta:=2\int_\beta c_1(T_X),
\qquad
\deg t_i:=2-d_i.
\]
The dimension axiom implies that $\star_{t'}$ is homogeneous of degree zero. More precisely, if
\[
T_i\star_{t'}T_j=\sum_k C_{ij}^{k}T_k,
\]
then every monomial occurring in $C_{ij}^{k}$ has weight $d_i+d_j-d_k$. The \emph{Euler vector field}
\[
\mathsf{Eu}_{t'}:=c_1(T_X)+\sum_{i\not\in\{1,\ldots,\rho\}} \left(1-\frac{1}{2}\deg T_i\right)t_i T_i \in QH^{\mathrm{even}}(X),
\]
upon writing $c_1(T_X)=\sum_{i=1}^\rho r_i T_i$ ($r_i\in\mathbb Z$), corresponds under the identification above to the logarithmic derivation
\[
\sum_{i=1}^\rho r_i q_i\partial_{q_i}+\sum_{i\not\in\{1,\ldots,\rho\}} \left(1-\frac{1}{2}\deg T_i\right)t_i\partial_{t_i}.
\]
Denoting by $\mathsf{Lie}$ the Lie derivative, a direct computation using the degree homogeneity of $\star_{t'}$ shows that
\[
\mathsf{Lie}_{\mathsf{Eu}_{t'}}(\star_{t'})=\star_{t'}.
\]

\medskip

Throughout, we work over the algebraic closure $\mathbb K_X:=\overline{\operatorname{Frac}(\mathsf{Nov}_X^{\mathrm{red}})}$ of the fraction field of $\mathsf{Nov}_X^{\mathrm{red}}$, abbreviate $H^{\mathrm{even}}(X,\mathbb K_X):=H^{\mathrm{even}}(X,\mathbb Q)\otimes_{\mathbb Q}\mathbb K_X$, and consider the $\mathbb K_X$-linear multiplication operator
\[
\boldsymbol\kappa_{t'}
:=
\bigl(\mathsf{Eu}_{t'}\star_{t'}(-)\bigr)\in\operatorname{End}_{\mathbb K_X}\bigl(H^{\mathrm{even}}(X,\mathbb K_X)\bigr).
\]
By the \emph{multiplicity} of an eigenvalue of $\boldsymbol\kappa_{t'}$ in $\mathbb K_X$, we always mean its algebraic multiplicity as a root of the characteristic polynomial of $\boldsymbol\kappa_{t'}$.
\begin{remark}
The construction above is the even, reduced formal version of the base of the $\mathsf A$-model F-bundle of \cite[Sections~3.5.1--3.5.2]{KKPY}. 
\end{remark}

The main rationality obstruction we shall use comes from the theory of Hodge atoms of Katzarkov--Kontsevich--Pantev--Yu in the following simplified form observed by Cai \cite{cai2026quarticthreefoldsymplecticallyirrational}:
\begin{proposition}[{\cite[Theorem~5.6]{cai2026quarticthreefoldsymplecticallyirrational}, after \cite{KKPY}}]\label{thm:cai-bound}
For every smooth projective rational threefold $X$, every eigenvalue of $\boldsymbol{\kappa}_{t'}$ has algebraic multiplicity at most two.
\end{proposition}

\begin{remark}
All completions used in this paper are ordinary completions. The proof of Theorem \ref{thm:cai-bound}, however, relies on the blowup formula of \cite{Iritani}, which is formulated using graded completions; see \cite[Remark~1.3 and Section~2.2]{Iritani}. This distinction causes no difficulty here. By the base-change comparison of Section~\ref{sec:coefficient-comparison}, the algebraic multiplicities of the eigenvalues of the Euler multiplication operator are preserved after passing between the corresponding coefficient fields.
\end{remark}

\subsection{Local Artinian factors}

We shall use the algebra structure of quantum cohomology to obtain lower bounds for the algebraic multiplicities of the eigenvalues of $\boldsymbol\kappa_{t'}$. The relevant observation is that multiplication by any element of a finite-dimensional local algebra over an algebraically closed field has a single eigenvalue.

Of interest to us is the finite-dimensional commutative unital $\mathbb K_X$-algebra
\[
A_X:=
QH^{\mathrm{even}}(X)
\otimes_{\mathsf{Nov}^{\mathrm{red}}_X}\mathbb K_X.
\]
As a finite-dimensional commutative algebra, $A_X$ is Artinian and has a decomposition, unique up to permutation,
\[
A_X\simeq\prod_{\nu=1}^{s}A_\nu,
\]
where the $A_\nu$ are local Artinian $\mathbb K_X$-algebras. We call them the \emph{local factors} of $A_X$. This is the standard decomposition obtained from the Chinese remainder theorem and nilpotence of the Jacobson radical.

Let $\mathfrak m_\nu$ be the maximal ideal of $A_\nu$.
Since $\mathbb K_X$ is algebraically closed,
\[
A_\nu/\mathfrak m_\nu\simeq\mathbb K_X.
\]
Write the image of $\mathsf{Eu}_{t'}$ in $A_\nu$ as
\[
\mathsf{Eu}_\nu=\lambda_\nu\mathbf1_\nu+n_\nu,
\qquad
\lambda_\nu\in\mathbb K_X,\quad n_\nu\in\mathfrak m_\nu,
\]
where $\mathbf1_\nu$ is the unit of $A_\nu$. Using the product decomposition to identify $A_\nu$ with a summand of $A_X$, we obtain
\[
\left.\boldsymbol\kappa_{t'}\right|_{A_\nu}
=\lambda_\nu\operatorname{id}_{A_\nu}+L_{n_\nu}.
\]
Here $L_{n_\nu}$ denotes multiplication by $n_\nu$ and is nilpotent because $\mathfrak m_\nu$ is nilpotent. Consequently,
\[
\det\left(
z\,\operatorname{id}_{A_X}-\boldsymbol\kappa_{t'}
\right)
=
\prod_{\nu=1}^{s}
(z-\lambda_\nu)^{\dim_{\mathbb K_X}A_\nu}.
\]

More precisely, the generalized eigenspace for an eigenvalue
$\lambda$ is
\[
\bigoplus_{\lambda_\nu=\lambda}A_\nu.
\]
Indeed, $\boldsymbol\kappa_{t'}-\lambda\operatorname{id}$ is nilpotent on a factor with $\lambda_\nu=\lambda$ and invertible on every other factor. Thus a local factor coincides with its corresponding generalized eigenspace if and only if its eigenvalue occurs on no other local factor.

In particular, a local factor of dimension $d$ gives an eigenvalue of algebraic multiplicity at least $d$. By Proposition~\ref{thm:cai-bound}, a local factor of dimension at
least three obstructs rationality for a smooth projective threefold.

\medskip

The arguments in Section~\ref{sec:conic-bundle-quantum} require changes of coefficient rings. We conclude this section by remarking that the injective changes considered below preserve the dimensions of local factors and their radical powers, as well as the eigenvalue multiplicities of $\boldsymbol\kappa_{t'}$, after passing to algebraically closed fraction fields.

\subsection{Changing coefficient rings}
\label{sec:coefficient-comparison}

Let
\[
\mathrm{Nov}_X^\circ
:=
\mathbb{Q}\bigl[C_{ij}^{k},(\boldsymbol{\kappa}_{t'})_{ij}\bigr]
\subset \mathrm{Nov}_X^{\mathrm{red}},
\]
where as above the entries $C_{ij}^{k}$ are the structure constants for quantum multiplication taken in the chosen homogeneous basis. The quantum algebra and its Euler multiplication operator are defined
over this subalgebra. Suppose that $\iota:\mathsf{Nov}^{\circ}_X\hookrightarrow\mathcal O$ is an injective homomorphism into an integral domain. Applying
$\iota$ to these coefficients and extending scalars to $\overline{\operatorname{Frac}(\mathcal O)}$ gives another algebra with a multiplication operator. The dimensions of its local Artinian factors and their radical powers agree with those of the corresponding factors over $\mathbb K_X$. The algebraic multiplicities of the corresponding eigenvalues also agree.

Indeed, after choosing compatible field embeddings, both constructions are obtained from the same algebra and operator
over $\overline{\operatorname{Frac}(\mathsf{Nov}^{\circ}_X)}$ by extending algebraically closed fields. Such extensions preserve the local Artinian decomposition, commute with taking radical powers, and preserve root multiplicities of characteristic polynomials.

\section{Quantum cohomology of standard conic bundles}
\label{sec:conic-bundle-quantum}

Throughout this section, let $\pi: X\to S$ be a projective standard conic bundle over a smooth projective rational surface $S$, with nonempty discriminant curve $\Delta\subset S$. Recall that this means that $X$ is smooth, $\pi$ is a flat proper morphism with plane-conic fibers, and $\rho(X/S)=1$. Set
\[
\Delta_S:=2K_S+\Delta,\qquad
\rho_S:=\operatorname{rank}\mathsf{NS}(S)=b_2(S).
\]
Then $\rho:=\operatorname{rank}\mathsf{NS}(X)=\rho_S+1$. In this section, we prove:
\begin{theorem}\label{thm:conic-bundle-local-factor}
Let $\pi: X\to S$ be a projective standard conic bundle over a smooth projective rational surface $S$. Assume that $\Delta_S=2K_S+\Delta$ is nef. Then
\[
QH^{\mathrm{even}}(X)
\otimes_{\mathsf{Nov}^{\mathrm{red}}_X}\mathbb K_X
\]
has a local Artinian factor $\mathcal H_-$ satisfying
\[
\begin{aligned}
\dim_{\mathbb K_X}\mathcal H_-&=\rho_S+2,&
\dim_{\mathbb K_X}\operatorname{rad}(\mathcal H_-)&=\rho_S+1,\\
\dim_{\mathbb K_X}\bigl(\operatorname{rad}(\mathcal H_-)^2\bigr)&=1,&
\operatorname{rad}(\mathcal H_-)^3&=0.
\end{aligned}
\]
Moreover, $\boldsymbol\kappa_{t'}$ has an eigenvalue of algebraic multiplicity exactly $\rho_S+2\geq3$, and so $X$ is irrational.
\end{theorem}

\subsection{Classical cohomology}

Since the fibers of $\pi$ are plane conics, we have
\[
\pi_*\mathcal O_X=\mathcal O_S,\qquad
R^j\pi_*\mathcal O_X=0\quad(j>0).
\]
As $S$ is rational, it follows that $H^2(X,\mathcal O_X)=0$.
The Lefschetz $(1,1)$-theorem, standardness, and Poincar\'e duality give
\[
b_2(X)=b_4(X)=\rho_S+1,\qquad
\dim_{\mathbb Q}H^{\mathrm{even}}(X,\mathbb Q)=2\rho_S+4.
\]

Let $\ell\in N_1(X,\mathbb Z)$ be the class of a line component of a
reducible fiber over a smooth point of $\Delta$. The relative
anticanonical embedding gives
\[
\int_\ell c_1(T_X)=1.
\]
Since $\rho(X/S)=1$, all such line components have the same rational
homology class, and a smooth fiber has class $2\ell$. Moreover,
\[
\ker\bigl(\pi_*:N_1(X,\mathbb Z)\to N_1(S,\mathbb Z)\bigr)
=\mathbb Z\ell.
\]
Indeed, the rational kernel is $\mathbb Q\ell$, and if $a\ell\in N_1(X,\mathbb Z)$, then $a=\int_{a\ell}c_1(T_X)\in\mathbb Z$.

Write
\[
H_S:=H^{\mathrm{even}}(S,\mathbb Q)
=\mathbb Q\mathbf1\oplus D_S\oplus\mathbb Qp,\qquad
D_S:=H^2(S,\mathbb Q),
\]
where $p=[\mathrm{pt}]$ is the point class normalized by $\int_Sp=1$. Fix a $\mathbb Q$-basis $\Delta_1,\ldots,\Delta_{\rho_S}$ of $D_S$, and put
\[
g_{ij}:=\int_S\Delta_i\smile\Delta_j.
\]
Then $(g_{ij})$ is nondegenerate, and
\[
\Delta_i\smile\Delta_j=g_{ij}p,\qquad
\Delta_i\smile p=p^{\smile2}=0.
\]
Hereafter, we use the same notation for a class on $S$ and its pullback to $X$. Moreover, when no ambiguity arises, we will use juxtaposition to denote the classical cup product.

Let
\[
\omega_{X/S}:=\omega_X\otimes\pi^*\omega_S^{-1},
\qquad
h:=c_1(\omega_{X/S}^{-1})
=c_1(T_X)-\pi^*c_1(T_S).
\]
Since $h$ has degree two on a smooth fiber, $\pi_*h=2$. By \cite[Lemma~3.11]{Prokhorov},
\[
\pi_*(K_X^2)\equiv-4K_S-\Delta.
\]
Because numerically equivalent divisors on $S$ have the same rational cohomology class, the former can be understood as an equality in $H^2(S,\mathbb Q)$. Since $h=-K_X+\pi^*K_S$ and $\pi_*K_X=-2$, the projection formula gives
\[
\pi_*(h^2)
=\pi_*(K_X^2)-2K_S\,\pi_*K_X
=-\Delta.
\]
Define
\[
x:=\frac12h+\frac18\Delta.
\]
The projection formula gives
\[
\pi_*x=1,\qquad
\pi_*(x^2)
=\frac14\pi_*(h^2)+\frac18\Delta\,\pi_*h=0.
\]
These identities give the following presentation.

\begin{lemma}\label{lem:cohomology-presentation}
There is an isomorphism of graded $H_S$-modules
\[
H^{\mathrm{even}}(X,\mathbb Q)=\pi^*H_S\oplus \left(x\smile\pi^*H_S\right).
\]
Moreover, there is a unique $\alpha\in\mathbb Q$ such that $x^2=\alpha p$. For $a,b,a_1,a_2,b_1,b_2\in H_S$, omitting pullbacks, we have
\[
\int_X(a+xb)=\int_Sb,\qquad
c_1(T_X)=2x+c_1(T_S)-\frac14\Delta,
\]
and
\[
(a_1+xb_1,a_2+xb_2)_X=\int_S(a_1b_2+b_1a_2).
\]
Here, integration of an inhomogeneous class means integration of its top-degree component.
\end{lemma}
\begin{proof}
Since $\pi_*(x\pi^*a)=a$, the pullback $\pi^*$ is injective. If $\pi^*a+x\pi^*b=0$, applying $\pi_*$ gives $b=0$, and then
$a=0$. The resulting injection
\[
\pi^*H_S\oplus x\pi^*H_S
\longrightarrow H^{\mathrm{even}}(X,\mathbb Q)
\]
is an isomorphism because both sides have dimension $2\rho_S+4$. In degree four, this identifies the kernel of $\pi_*$ with $\mathbb Qp$. Thus, $\pi_*(x^2)=0$ implies $x^2=\alpha p$ for a unique $\alpha\in\mathbb Q$.

The integration and pairing formulas follow from $\pi_*x=1$ and $\pi_*(x^2)=0$. The expression for $c_1(T_X)$ follows from the definition of $x$.
\end{proof}

We now choose coordinates adapted to this presentation. Since $c_1$ identifies $\mathsf{NS}(X)\otimes\mathbb Q$ with $H^2(X,\mathbb Q)$, choose the basis elements outside degree two to be
\[
\mathbf1,\quad p,\quad
x\Delta_1,\ldots,x\Delta_{\rho_S},\quad xp.
\]
Let $u,w,c^1,\ldots,c^{\rho_S},\xi$ be the corresponding formal coordinates. Then
\begin{equation}\label{eq:reduced-parameter}
t'=a'+xb',
\end{equation}
where
\[
a':=u\mathbf1+wp,\qquad
b':=c+\xi p,\qquad
c:=\sum_{i=1}^{\rho_S}c^i\Delta_i.
\]
The classes $\mathbf1,p,x\Delta_i,xp$ have real cohomological degrees $0,4,4,6$, respectively. Hence, the formal weights of their coordinates are
\[
\deg u=2,\qquad
\deg w=\deg c^i=-2,\qquad
\deg\xi=-4.
\]
Moreover,
\begin{equation}\label{eq:conic-euler}
\mathsf{Eu}_{t'}
=2x+c_1(T_S)-\frac14\Delta
+u\mathbf1-wp-xc-2\xi xp.
\end{equation}

\subsection{The vertical quantum product}
\label{sec:vertical}
A class $\beta\in\mathsf{NE}(X,\mathbb Z)$ is called \emph{vertical} if $\pi_*\beta=0$, equivalently if $\beta=m\ell$ for some $m\geq0$. Put $q_f:=q^\ell$. The \emph{vertical quantum product} $\star_{\mathrm{vert}}$ is obtained from $\star_{t'}$ by setting $q^\beta=0$ for all nonvertical classes. We should note that this is a ring specialization: a sum of effective classes is vertical only if both summands are vertical.

We compute this product and show that, after changing coefficients, its algebra splits into two copies of $H_S$. All products in $H_S$ below are cup products. Since $(b')^3=0$,
\[
e^{b'}=\mathbf1+b'+\frac12(b')^2,
\qquad
e^{\pm b'/2}=\mathbf1\pm\frac12b'+\frac18(b')^2.
\]

\begin{proposition}\label{prop:vertical-splitting}
\begin{enumerate}
\item[\textup{(i)}]
By multilinearity and the fundamental-class and divisor axioms, the positive-degree vertical primary genus-zero invariants with even insertions are determined by
\begin{equation}\label{eq:vertical-invariants}
\begin{aligned}
\langle xD\rangle_{0,1,\ell}&=\int_S\Delta D,\\
\langle xp\rangle_{0,1,2\ell}&=1,\\
\langle xD_1,xD_2\rangle_{0,2,2\ell}&=\int_SD_1D_2,
\end{aligned}
\qquad D,D_1,D_2\in D_S.
\end{equation}
Among invariants with insertions in $\{p,x\Delta_i,xp\}$, all other positive-degree vertical invariants vanish. For $a_1,a_2,b_1,b_2\in H_S$, it follows that
\[
(a_1+xb_1)\star_{\mathrm{vert}}(a_2+xb_2)
=a_1a_2+\mathcal B(b')b_1b_2+x(a_1b_2+a_2b_1),
\]
where
\begin{equation*}
\mathcal B(b'):=x\star_{\mathrm{vert}}x
=\alpha p+q_f^2e^{b'}+\frac14q_f\Delta e^{b'/2}.
\end{equation*}

\item[\textup{(ii)}]
The product in \textup{(i)} has coefficients in
\[
R_0:=\mathbb Q[q_f,u,w,c^1,\ldots,c^{\rho_S},\xi].
\]
Regard it as an algebra over $R_0$ and extend scalars to
\[
R_f:=k_f[[u,w,c^1,\ldots,c^{\rho_S},\xi]],
\qquad k_f:=\mathbb Q(q_f).
\]
Set
\[
\eta:=\alpha p-\frac1{64}\Delta^2=\gamma p,
\qquad
\gamma:=\alpha-\frac1{64}\int_S\Delta^2,
\]
and define
\begin{equation*}
\begin{aligned}
r(b')&:=q_fe^{b'/2}+\frac18\Delta
+\frac{\eta}{2q_f}e^{-b'/2},\\
\mathcal R(b')&:=2q_fe^{b'/2}+\frac18\Delta b'
-\frac{\eta}{q_f}e^{-b'/2},\\
t_\pm&:=a'\pm\mathcal R(b').
\end{aligned}
\end{equation*}
Then $r(b')^2=\mathcal B(b')$, and the map
\[
\Phi(a+xb):=
\bigl(a+r(b')b,\ a-r(b')b\bigr)
\]
is an $R_f$-algebra isomorphism from the vertical algebra to
\[
(H_S\otimes_{\mathbb Q}R_f)\times
(H_S\otimes_{\mathbb Q}R_f),
\]
where each factor is taken with the classical cup product.

\item[\textup{(iii)}]
Let $(\mathsf{Eu}_+,\mathsf{Eu}_-):=\Phi(\mathsf{Eu}_{t'})$.
Writing
\[
t_\pm=(t_\pm)_0\mathbf1+(t_\pm)_1+(t_\pm)_2p,
\qquad
(t_\pm)_1\in D_S\otimes_{\mathbb Q}R_f,
\]
we have
\begin{equation}\label{eq:Euler-pm}
\begin{aligned}
\mathsf{Eu}_+&=(t_+)_0\mathbf1+c_1(T_S)-(t_+)_2p,\\
\mathsf{Eu}_-&=(t_-)_0\mathbf1-\frac12\Delta_S-(t_-)_2p.
\end{aligned}
\end{equation}
\end{enumerate}
\end{proposition}
\begin{proof}
\textup{(i)}
A stable map of class $m\ell$, with $m>0$, is contained in a fiber of $\pi$. Let
\[
b_{m,n}:\overline{\mathcal M}_{0,n}(X,m\ell)\longrightarrow S
\]
record that fiber.

After applying the fundamental-class and divisor axioms, the remaining insertions are $p, xD,xp$, with $D\in D_S$. Write them as $\gamma_i=x^{\varepsilon_i}\pi^*\alpha_i$, where $\varepsilon_i\in\{0,1\}$. Since $\pi\circ\operatorname{ev}_i=b_{m,n}$, the product of evaluation classes contains $b_{m,n}^*(\alpha_1\cdots\alpha_n)$. If an invariant is nonzero, then $\sum_i\deg\alpha_i\leq4$.
The pairs $(\deg\gamma_i,\deg\alpha_i)$ are $(4,4),(4,2),(6,4)$, respectively, so the dimension axiom gives
\[
2m=\sum_i(\deg\gamma_i-2)
\leq\sum_i\deg\alpha_i\leq4.
\]
Thus, the only possibilities are the three types in \eqref{eq:vertical-invariants} and $\langle p\rangle_{0,1,\ell}$. If no insertions remain, the dimension axiom gives vanishing because $m>0$.

We first compute the class-$\ell$ invariants. Over the smooth locus of $\Delta$, such a map parametrizes a line component $L$ of a reducible fiber. Restricting $\pi$ to a curve transverse to $\Delta$ exhibits $L$ as a $(-1)$-curve in a smooth surface. The normal-bundle sequence splits, giving
\[
N_{L/X}\simeq
\mathcal O_{\mathbb P^1}(-1)\oplus\mathcal O_{\mathbb P^1}.
\]
Hence these maps are unobstructed, and their one-point evaluation map has generic degree one onto
$\pi^{-1}(\Delta)$. Contributions over $\operatorname{Sing}(\Delta)$ cannot change the pushforward two-cycle, since $\pi^{-1}(\operatorname{Sing}(\Delta))$ has dimension at most one. Thus
\[
\operatorname{ev}_*
[\overline{\mathcal M}_{0,1}(X,\ell)]^{\mathrm{vir}}
=[\pi^{-1}(\Delta)].
\]
Therefore
\[
\langle xD\rangle_{0,1,\ell}
=\int_XxD\Delta=\int_S\Delta D,
\qquad
\langle p\rangle_{0,1,\ell}
=\int_Xp\Delta=0.
\]

Over $S\setminus\Delta$, a map of class $2\ell$ has degree one onto a smooth fiber $F\simeq\mathbb P^1$, with trivial normal bundle and $\int_Fx=1$. For $j=1,2$, the fiberwise evaluation map to $F^j$ has generic degree one, so the intersection number of the $j$ evaluation classes of $x$ is $(\int_Fx)^j=1$. Thus
\[
(b_{2,j})_*\left(
\prod_{i=1}^{j}\operatorname{ev}_i^*x
\cap[\overline{\mathcal M}_{0,j}(X,2\ell)]^{\mathrm{vir}}
\right)=[S],
\qquad j=1,2.
\]
The equality holds globally because no two-dimensional cycle can be supported on $\Delta$. Pairing with $p$ for $j=1$ and with $D_1D_2$ for $j=2$ proves the remaining formulas in \eqref{eq:vertical-invariants}.

Divisors pulled back from $S$ have zero vertical degree, so their quantum multiplication is classical. Since divisors and the unit generate $H_S$, associativity gives
\[
a\star_{\mathrm{vert}}\gamma=a\gamma
\qquad
(a\in H_S,\ \gamma\in H^{\mathrm{even}}(X,\mathbb Q)).
\]
It remains to determine $\mathcal B:=x\star_{\mathrm{vert}}x$. For every $a\in H_S$, invariance of the pairing gives
\[
(\mathcal B,a)_X=(x,xa)_X=\int_Xx^2a=0.
\]
By Lemma~\ref{lem:cohomology-presentation}, this means that the coefficients of $\mathcal B$ lie in $H_S$.

Using the computed invariants and
$\int_\ell x=1/2$, $\int_{2\ell}x=1$, the divisor axiom yields
\[
\begin{aligned}
(\mathcal B,xp)_X
&=q_f^2,\\
(\mathcal B, xD)_X
&=q_f^2\int_ScD+\frac14q_f\int_S\Delta D,\\
(\mathcal B,x)_X
&=\alpha+q_f^2\xi
+\frac{q_f^2}{2}\int_Sc^2
+\frac{q_f}{8}\int_S\Delta c.
\end{aligned}
\]
These pairings determine
\[
\mathcal B
=\alpha p+q_f^2e^{b'}+\frac14q_f\Delta e^{b'/2}.
\]
The general product formula follows from the classical
$H_S$-action.

\textup{(ii)}
The structure constants in \textup{(i)} belong to $R_0$, so we may extend scalars to $R_f$. Since $\Delta\eta=\eta^2=0$,
\[
r(b')^2
=q_f^2e^{b'}+\frac14q_f\Delta e^{b'/2}
+\frac1{64}\Delta^2+\eta
=\mathcal B(b').
\]
Moreover, $r(b')$ is invertible: its degree-zero component is $q_f$, and its positive-degree part is nilpotent.

The vertical algebra is generated over $H_S\otimes R_f$ by $x$, with relation $x\star_{\mathrm{vert}}x=r(b')^2$. Evaluation at $x=\pm r(b')$ therefore defines $\Phi$. Its inverse is
\[
(A_+,A_-)\longmapsto
\frac{A_++A_-}{2}
+x\,\frac{A_+-A_-}{2r(b')},
\]
which proves the splitting.

\textup{(iii)}
Put
\[
\mathcal R_2:=
q_f\xi+\frac{q_f}{4}\int_Sc^2
+\frac18\int_S\Delta c-\frac{\gamma}{q_f}.
\]
Expanding the finite exponentials and using $D_Sp=p^2=0$ gives
\[
\mathcal R(b')
=2q_f\mathbf1+q_fc+\mathcal R_2p,
\]
and
\[
r(b')(2\mathbf1-c-2\xi p)
=2q_f\mathbf1+\frac14\Delta-\mathcal R_2p.
\]
Applying $\Phi$ to \eqref{eq:conic-euler}, we obtain
\[
\mathsf{Eu}_\pm
=(u\pm2q_f)\mathbf1
+\left(c_1(T_S)-\frac14\Delta\pm\frac14\Delta\right)
-(w\pm\mathcal R_2)p.
\]
Since
\[
(t_\pm)_0=u\pm2q_f,\qquad
(t_\pm)_2=w\pm\mathcal R_2,\qquad
c_1(T_S)-\frac12\Delta=-\frac12\Delta_S,
\]
this proves \eqref{eq:Euler-pm}.
\end{proof}

\subsection{Coefficients and divisor variables}

Our next step is to compare the vertical product with the full quantum product. We first choose coordinates in which the origin of the negative factor is at zero. We then choose a coefficient ring in which this change of coordinates is well-defined and the vertical product is recovered by setting the contributions of nonvertical curves equal to zero.

At $c=\xi=0$, Proposition~\ref{prop:vertical-splitting} gives
\[
t_-=(u-2q_f)\mathbf1+(w+\gamma/q_f)p.
\]
Thus the origin of this factor corresponds to $u=2q_f$ and $w=-\gamma/q_f$. We introduce coordinates centered
at this point:
\begin{equation}\label{eq:uw-translation}
u=2q_f+u',\qquad
w=-\frac{\gamma}{q_f}+w'.
\end{equation}
In particular, $t_-=0$ when $u'=w'=c=\xi=0$.

This substitution requires some care. A formal power series in $w$ cannot, in general, be evaluated at $w=-\gamma/q_f+w'$, because its constant coefficient in $w'$ could involve an undefined infinite sum. We shall show that the coefficients occurring in our quantum product are polynomial in $w$ after fixing the projected curve class (making the substitution legitimate).

\medskip

To keep track of projected curve classes, put
\[
L:=\pi_*N_1(X,\mathbb Z),\qquad
P:=\pi_*\mathsf{NE}(X,\mathbb Z)\subseteq L.
\]
Thus $L$ is the group of curve classes on $S$ obtained by pushforward, and $P$ consists of the classes obtained from effective curves on $X$.

Since
\[
\ker\bigl(\pi_*:N_1(X,\mathbb Z)\to L\bigr)=\mathbb Z\ell
\]
and $L$ is a free abelian group, we may choose a homomorphism
\[
\sigma:L\longrightarrow N_1(X,\mathbb Z),
\qquad \pi_*\circ\sigma=\operatorname{id}_L.
\]
For example, choose a basis of $L$ and lift each basis element
to $N_1(X,\mathbb Z)$. Every curve class $B$ then has a unique
expression
\[
B=\sigma(\beta)+m\ell,\qquad
\beta=\pi_*B,\quad m\in\mathbb Z.
\]
We accordingly write
\[
q^B=q_S^\beta q_f^m.
\]
The symbols $q_S^\beta$ record the projected class. The chosen lift $\sigma(\beta)$ need not be effective, so $m$ may be negative even when $B$ is effective. Fix a very ample divisor $A_S$ and set
\[
l(\beta):=\int_\beta A_S.
\]
Every nonzero class in $P$ has positive $l$-degree, and only finitely many classes in $P$ have bounded $l$-degree.

\medskip

We now justify the required polynomial dependence on $w$. Recall that $w$ is the parameter multiplying the pulled-back point class $p$. 

Fix $\beta\in P\setminus\{0\}$. The pushforward of a stable map with projected class $\beta$ is an effective one-cycle of degree $l(\beta)$ on $S$. Since the domain is connected and $\beta\neq0$, all projected marked points lie on the support of this cycle, including those carried by components contracted by $\pi$. Let
\[
M_\beta:=
\dim\operatorname{Chow}_{1,l(\beta)}(S),
\]
where the Chow variety parametrizes effective one-cycles of this degree for the embedding defined by $A_S$.

Consider the incidence variety parametrizing such a cycle together with $k$ points on its support. Its dimension is at most $M_\beta+k$, whereas $\dim S^k=2k$. If $k>M_\beta$, its image cannot fill $S^k$. Consequently, no such cycle passes through $k$ general prescribed points of $S$. It follows that every Gromov--Witten invariant of projected class $\beta$ with more than $M_\beta$ insertions of $p$ vanishes.

Thus, at fixed $\beta\neq0$, the structure constants of the quantum product are polynomial in $w$. For $\beta=0$, this follows from the explicit vertical calculation. The fundamental-class axiom also shows that the product is independent of $u$. Since the Euler vector field is affine in $u,w$, the matrix entries of $\boldsymbol\kappa_{t'}$ are polynomial in these two variables as well, at each fixed projected class.

There is likewise no infinite summation over the fiber degree once the projected class and parameter monomial are fixed. Indeed, homogeneity determines the fiber degree. If a homogeneous coefficient has degree $d$, then a monomial
\[
q_S^\beta q_f^m u^a w^b
(c^1)^{I_1}\cdots(c^{\rho_S})^{I_{\rho_S}}\xi^j
\]
can occur only if
\[
d=
2\int_{\sigma(\beta)}c_1(T_X)
+2m+2a-2b-2|I|-4j,
\qquad |I|:=\sum_i I_i.
\]
Thus, for fixed $\beta$ and parameter monomial, the exponent $m$ is uniquely determined.

Recall that $\mathsf{Nov}^{\circ}_X$ is the $\mathbb Q$-algebra generated by these structure constants and matrix entries. The preceding observations imply that, for every $f\in\mathsf{Nov}^{\circ}_X$, we can write
\[
f=\sum_{\beta\in P}q_S^\beta f_\beta,
\qquad
f_\beta\in
\mathbb Q[q_f,q_f^{-1}]
[[c^1,\ldots,c^{\rho_S},\xi]][u,w].
\]
The brackets $[u,w]$ mean that $f_\beta$ is a polynomial in $u,w$, although its coefficients may be formal power series in $c^1,\ldots,c^{\rho_S},\xi$. This property is preserved under products because each $\beta\in P$ has only finitely many decompositions into effective projected classes.

\medskip

We can now specify the new coefficient ring. Treat $q_f$ as an invertible coefficient and put
\[
k_f:=\mathbb Q(q_f),\qquad
R:=k_f[[P]]_l
=
\left\{
\sum_{\beta\in P}a_\beta q_S^\beta:
a_\beta\in k_f
\right\}.
\]
Multiplication is given by
$q_S^\beta q_S^{\beta'}=q_S^{\beta+\beta'}$. We complete with respect to the degree filtration
\[
F^rR:=
\left\{
\sum_\beta a_\beta q_S^\beta:
a_\beta=0\text{ whenever }l(\beta)<r
\right\},
\qquad r\geq0.
\]
Thus, $F^rR$ contains contributions with projected degree at least $r$. This makes $R$ a complete integral domain, and
\[
R/F^1R=k_f,\qquad
F^rR\,F^sR\subseteq F^{r+s}R.
\]
In particular, reduction modulo $F^1R$ discards all nonvertical curve classes.

Writing $\mathbf c=(c^1,\ldots,c^{\rho_S})$, we define
\[
\mathsf{Nov}^{\circ}_X
\longrightarrow R[[u',w',\mathbf c,\xi]]
\]
by applying \eqref{eq:uw-translation} to each polynomial $f_\beta$. At each fixed $\beta$ and parameter monomial, the substitution involves only finitely many terms, so this is well-defined. It is injective because polynomial translation has the inverse
\[
u'=u-2q_f,\qquad w'=w+\gamma/q_f.
\]
We have therefore justified the change of origin on the algebra containing all coefficients needed for the quantum product and $\boldsymbol\kappa_{t'}$.

\medskip

We next restore the divisor parameters. This allows us to represent all cohomology classes by coordinate vector fields, including the divisor classes that were absorbed into the Novikov variables in the reduced description.

Introduce formal coordinates
$v^1,\ldots,v^{\rho_S},z$ corresponding to
$\Delta_1,\ldots,\Delta_{\rho_S},x$, and put
\[
v:=\sum_i v^i\Delta_i,\qquad
\mathbf v:=(v^1,\ldots,v^{\rho_S}),\qquad
\widehat{\mathcal O}:=
R[[u',w',\mathbf v,z,\mathbf c,\xi]].
\]
The full formal parameter is
\begin{equation}\label{eq:full-conic-parameter}
\begin{aligned}
t&=a+xb,\\
a&=(2q_f+u')\mathbf1+v+
   \left(-\frac{\gamma}{q_f}+w'\right)p,\\
b&=z\mathbf1+c+\xi p.
\end{aligned}
\end{equation}

The divisor axiom describes exactly how the new parameters enter the product. In curve class $B$, inserting $v+zx$ multiplies an invariant by
\[
\int_B(v+zx)=\int_{\pi_*B}v+z\int_Bx.
\]
Summing over any number of these insertions gives the
exponential of this expression. Consequently, the coefficients
of the quantum product at $t$ are obtained by the simultaneous
substitution
\begin{equation}\label{eq:conic-coefficient-embedding}
\begin{aligned}
q^B&\longmapsto
q_S^\beta q_f^m
\exp\left(\int_\beta v+z\int_Bx\right),
\qquad B=\sigma(\beta)+m\ell,\\
u&\longmapsto2q_f+u',\\
w&\longmapsto-\frac{\gamma}{q_f}+w',
\end{aligned}
\end{equation}
with $c^i,\xi$ unchanged.

The preceding finiteness arguments make this substitution well-defined on $\mathsf{Nov}^{\circ}_X$. Setting $\mathbf v=z=0$ recovers the injective translation already constructed, so \eqref{eq:conic-coefficient-embedding} is injective as well.

 Since divisor parameters have coefficient zero in the linear part of the Euler vector field, the same substitution sends
the matrix of $\boldsymbol\kappa_{t'}$ to the matrix of multiplication by the Euler vector field at $t$.
Here $q_f$ is a coefficient in $R$, independent of $z$: the chosen center remains $u=2q_f$, $w=-\gamma/q_f$.

Finally, we specify the topology used for subsequent formal coordinate changes. For a finite tuple of variables $\mathbf y$, let $F^rR[[\mathbf y]]$ consist of series whose
coefficients all belong to $F^rR$. On $R[[\mathbf y]]$ we use the neighborhoods
\[
F^rR[[\mathbf y]]+(\mathbf y)^r,\qquad r\geq1.
\]
Thus terms tend to zero when either their projected curve
degree or their total parameter degree tends to infinity. In particular, a translation by a coefficient in $F^1R$ is well-defined, and successive coordinate corrections of increasing $F$-order converge in this topology.

\subsection{Extending the vertical decomposition}

We now work over the coefficient ring $\widehat{\mathcal O}$ constructed above. Reduction modulo $F^1$ retains exactly the vertical curve classes. Our first goal is to extend the two summands of the vertical quantum algebra to the full quantum product. We shall then choose coordinates on these summands and study the negative one.

Write $\star$ for the product at
\eqref{eq:full-conic-parameter}. The restored divisor coordinates allow us to regard this product as a multiplication of coordinate vector fields. More precisely, let
\[
\mathcal T:=\operatorname{Der}^{\mathrm{cont}}_R
(\widehat{\mathcal O})
\simeq H^{\mathrm{even}}(X,\mathbb Q)
\otimes_{\mathbb Q}\widehat{\mathcal O},
\]
where an $R$-linear derivation fixes all coefficients in $R$. The identification is
\[
\begin{gathered}
\partial_{u'}\longleftrightarrow\mathbf1,\qquad
\partial_{v^i}\longleftrightarrow\Delta_i,\qquad
\partial_z\longleftrightarrow x,\\
\partial_{w'}\longleftrightarrow p,\qquad
\partial_{c^i}\longleftrightarrow x\Delta_i,\qquad
\partial_\xi\longleftrightarrow xp.
\end{gathered}
\]
In particular, the unit of $\star$ is $\partial_{u'}$, which we also denote by $\mathbf1$ under this identification.

\medskip

We first describe the vertical decomposition in these coordinates. The formulas of Proposition~\ref{prop:vertical-splitting} extend
to $b=z\mathbf1+c+\xi p$ by putting
\[
\begin{aligned}
r(b)&:=q_fe^{b/2}+\frac18\Delta
       +\frac{\eta}{2q_f}e^{-b/2},\\
\mathcal R(b)&:=2q_fe^{b/2}+\frac18\Delta b
       -\frac{\eta}{q_f}e^{-b/2},\\
t_\pm&:=a\pm\mathcal R(b).
\end{aligned}
\]
Here $\eta=\gamma p$, as before, and
$e^{\pm b/2}=e^{\pm z/2}e^{\pm(c+\xi p)/2}$.
The scalar component of $r(b)$ is $q_fe^{z/2}$, so $r(b)$ is
invertible. Thus the vertical algebra has the two idempotents
\[
\epsilon_\pm^{(0)}
=\frac12\bigl(\mathbf1\pm r(b)^{-1}x\bigr).
\]
They satisfy
\[
\epsilon_\pm^{(0)}\star_{\mathrm{vert}}\epsilon_\pm^{(0)}
=\epsilon_\pm^{(0)},\qquad
\epsilon_+^{(0)}\star_{\mathrm{vert}}\epsilon_-^{(0)}=0,
\qquad \epsilon_+^{(0)}+\epsilon_-^{(0)}=\mathbf1.
\]
Multiplication by these idempotents projects onto the two vertical summands.

The same decomposition can be expressed in terms of a change of coordinates. Indeed,
\[
d\mathcal R_b(\delta b)=r(b)\delta b,
\]
so the differential of $(a,b)\mapsto(t_+,t_-)$ is
\[
(\delta a,\delta b)\longmapsto
\bigl(\delta a+r(b)\delta b,
      \delta a-r(b)\delta b\bigr).
\]
This is precisely the vertical algebra isomorphism of Proposition~\ref{prop:vertical-splitting}. In these coordinates, each summand has the constant cup product of $H_S$. At the chosen origin,
\[
t_-=0,\qquad
t_+=4q_f\mathbf1-\frac{2\gamma}{q_f}p.
\]
Hence $t_-$ and $t_+-(4q_f\mathbf1-2\gamma p/q_f)$ are coordinates centered at zero for the two vertical factors.

\medskip

We now restore the contributions of nonvertical curves. The idempotents $\epsilon_\pm^{(0)}$ lift uniquely to idempotents $\epsilon_\pm$ for $\star$. To see this, solve
\[
\epsilon\star\epsilon-\epsilon=0
\]
successively in the filtration $F^r$. The linear equation for a correction $h$ has the form
\[
(2\epsilon-\mathbf1)\star h
=-(\epsilon\star\epsilon-\epsilon).
\]
Its coefficient is invertible: modulo $F^1$,
$(2\epsilon_\pm^{(0)}-\mathbf1)^2=\mathbf1$.
Each correction raises the order of the error, and completeness gives the desired lift. Uniqueness also gives
$\epsilon_- =\mathbf1-\epsilon_+$, and therefore
\[
\epsilon_+\star\epsilon_-=0,
\qquad
\mathcal T=\mathcal T_+\oplus\mathcal T_-,
\qquad
\mathcal T_\pm:=\epsilon_\pm\star\mathcal T.
\]
Both summands have rank $\rho_S+2$, and products between them vanish.

This gives a decomposition of the algebra of vector fields. To achieve this using two groups of coordinates, we use the compatibility between the quantum product and the Lie bracket. For vector fields $Y,V,W$, define
\[
\mathcal P_Y(V,W):=
[Y,V\star W]-[Y,V]\star W-V\star[Y,W].
\]
Thus $\mathcal P_Y=\mathsf{Lie}_Y(\star)$ measures the change of the product along $Y$. The quantum product satisfies the Hertling--Manin identity
\begin{equation}\label{eq:conic-HM-identity}
\mathcal P_{Y\star Z}
=Y\star\mathcal P_Z+Z\star\mathcal P_Y;
\end{equation}
see \cite[Theorem~2(a)]{HertlingManin}. A formal power-series space over $R$, equipped with a commutative, associative, unital tangent multiplication satisfying this identity, is called a \emph{formal $F$-manifold over $R$}.

For an idempotent $\epsilon$, the identity gives $\mathcal P_\epsilon=2\epsilon\star\mathcal P_\epsilon$. Since multiplication by $\mathbf1-2\epsilon$ is invertible,
$\mathcal P_\epsilon=0$. If $Y,Z\in\mathcal T_+$, applying the identity to $\epsilon_-\star Y=0$ and evaluating on
$(\epsilon_-,Z)$ gives
\[
0=\epsilon_-\star\mathcal P_Y(\epsilon_-,Z)
=-\epsilon_-\star[Y,Z].
\]
Hence $[Y,Z]\in\mathcal T_+$. The same argument applies to
$\mathcal T_-$, so both summands are closed under Lie brackets. The formal Frobenius theorem therefore supplies coordinates in which disjoint groups of coordinate vector fields span these two summands.

Moreover, the product within either group is independent of the other group's coordinates. For example, let $A, B$ be coordinate fields in $\mathcal T_+$ and $Z$ a coordinate field in $\mathcal T_-$. Since $Z=\epsilon_-\star Z$ and $\mathcal P_{\epsilon_-}=0$, we have
\[
\mathcal P_Z(A,B)=\epsilon_-\star\mathcal P_Z(A,B).
\]
The coordinate fields commute, so the left-hand side is
$[Z,A\star B]$, which belongs to $\mathcal T_+$. The right-hand side belongs to $\mathcal T_-$. Both sides therefore vanish.

This is the decomposition argument in the proof of \cite[Theorem~4.2]{Hertling1999}, applied formally over $R$. The coordinate equations can be solved compatibly modulo $F^r$ by formal integration; the required positive integers are invertible because the coefficient field has characteristic zero. Passing to the limit gives coordinates that reduce to the vertical coordinates chosen above. We have thus obtained two formal $F$-manifolds whose product is the full formal space.

\medskip

We next construct a derivation suited to the negative factor. The reason for modifying the Euler vector field is already visible in the vertical calculation: its negative component contains the constant divisor term \[
C:=-\frac12\Delta_S=c_1(T_S)-\frac12\Delta.
\]
We shall replace differentiation in this divisor direction by an equivalent action on the curve coefficients. Nefness will then control that action.

For $D=\sum_iD^i\Delta_i\in D_S$, put
\[
\partial_{v,D}:=\sum_iD^i\partial_{v^i},
\qquad
\delta_D(q_S^\beta):=(D\cdot\beta)q_S^\beta.
\]
Here $\delta_D$ is a $k_f$-linear derivation of $R$; extend it
to $\widehat{\mathcal O}$ by fixing all formal coordinates.
The divisor substitution
\eqref{eq:conic-coefficient-embedding} shows that differentiating
a product coefficient along $\partial_{v, D}$ has the
same effect as applying $\delta_D$. Hence
\[
\mathsf{Lie}_{\partial_{v,D}-\delta_D}(\star)=0.
\]

Define
\[
\mathcal S_C:=\partial_{v,C}-\delta_C,
\qquad
\mathcal G:=\mathsf{Eu}_t-\mathcal S_C.
\]
In our coordinates, the Euler vector field is
\[
\mathsf{Eu}_t
=u\partial_{u'}+
\partial_{v,c_1(T_S)-\Delta/4}+2\partial_z
-w\partial_{w'}-\sum_i c^i\partial_{c^i}
-2\xi\partial_\xi,
\]
where $u=2q_f+u'$ and $w=-\gamma/q_f+w'$. It fixes $R$. Consequently,
\begin{equation}\label{eq:G-conformal}
\mathcal G|_R=\delta_C,
\qquad
\mathsf{Lie}_{\mathcal G}(\star)=\star.
\end{equation}
Although $\mathcal G$ need not be $R$-linear, its bracket with an $R$-linear vector field is $R$-linear. Its Lie derivative on the tangent multiplication is therefore well-defined.

The derivation $\mathcal G$ preserves both summands. Indeed, applying \eqref{eq:G-conformal} to an idempotent $\epsilon$ gives
\[
[\mathcal G,\epsilon]
-2\epsilon\star[\mathcal G,\epsilon]=\epsilon,
\]
and hence $[\mathcal G,\epsilon]=-\epsilon$. For any $V\in\mathcal T$, it follows that
\[
[\mathcal G,\epsilon\star V]
=\epsilon\star[\mathcal G,V].
\]
In product coordinates, this means that the coefficients of $\mathcal G$ in either group's coordinate directions depend only on that group's coordinates and on $R$. Thus $\mathcal G$ induces a derivation on each formal factor, acting on $R$ through $\delta_C$.

To compute the induced derivation on the negative factor modulo $F^1$, write
\[
t_-=U\mathbf1+\sum_iV_i\Delta_i+Wp,
\qquad C=\sum_iC^i\Delta_i.
\]
The calculation of \eqref{eq:Euler-pm} remains valid with $b=z\mathbf1+c+\xi p$. It gives the negative component of the
Euler vector field as
\[
U\partial_U+\sum_iC^i\partial_{V_i}-W\partial_W.
\]
Indeed, replacing $q_f$ by $q_fe^{z/2}$ gives the same scalar and degree-four components; the additional term $z\Delta/8$ in $\mathcal R(b)$ affects only the degree-two component of $t_-$.

Since $\partial_{v,C}$ translates $t_-$ by $C$, subtracting it removes $\sum_iC^i\partial_{V_i}$. The remaining term $\delta_C$ acts on curve coefficients and fixes the vertical coordinates. Choosing coordinates on the full negative factor that extend $U, V_i, W$, we obtain its coordinate ring
\[
\mathcal O_-:=R[[U,V_1,\ldots,V_{\rho_S},W]].
\]
Write $\circ$ for its tangent product and $\mathcal G_-$ for the induced derivation. Extending this tangent algebra from $\mathcal O_-$ to the full coordinate ring recovers $\mathcal T_-$.
We have
\begin{equation}\label{eq:negative-leading-field}
\mathcal G_-\equiv\mathcal G_0
:=\delta_C+U\partial_U-W\partial_W\pmod{F^1},
\qquad
\mathsf{Lie}_{\mathcal G_-}(\circ)=\circ.
\end{equation}
The congruence concerns coefficients in the coordinate directions; both derivations restrict to $\delta_C$ on $R$. Modulo $F^1$, the product is the cup product of $H_S$:
\[
\partial_{V_i}\circ\partial_{V_j}=g_{ij}\partial_W,
\qquad
\partial_{V_i}\circ\partial_W
=\partial_W\circ\partial_W=0,
\]
and its unit reduces to $\partial_U$.

\subsection{Nefness and non-splitting}

The preceding construction gives two summands of the full quantum algebra. We must still prove that the negative summand is local after passing to the algebraically closed fraction field. 

\begin{remark}
For the discussion that follows, it is important to note that locality of the negative factor in vertical quantum cohomology  does not immediately imply locality after the contributions of the remaining curve classes are restored. A simple toy model is
$\mathbb Q[[\tau]][\zeta]/(\zeta^2-\tau)$, which specializes at $\tau=0$ to a local algebra, but splits into two factors over an algebraic closure of $\mathbb Q((\tau))$.
\end{remark}

We will see that nefness rules out the above kind of splitting in our situation. Namely, $\delta_C$ acts on a projected curve monomial by
\[
\delta_C(q_S^\beta)=\lambda_\beta q_S^\beta,
\qquad
\lambda_\beta:=C\cdot\beta
=-\frac12\Delta_S\cdot\beta\leq0.
\]

The next lemma combines this inequality with \eqref{eq:negative-leading-field} to constrain the possible nonvertical corrections to the product. It shows that, after a formal change of coordinates, the products among the nonunit directions retain the same form as in $H_S$, although their coefficients may change.

\begin{lemma}\label{lem:nef-nonsplitting}
Let $R=k_f[[P]]_l$, and let $\delta$ be the $k_f$-linear
derivation determined by
\[
\delta(q_S^\beta)=\lambda_\beta q_S^\beta,
\]
where $\lambda:P\to\mathbb Q$ is additive and $\lambda_\beta\leq0$. Let
\[
\mathcal O:=R[[U,V_1,\ldots,V_s,W]],\qquad s\geq1,
\]
carry a formal $F$-manifold product $\circ$ with unit $e$. Suppose that modulo $F^1$ its product is
\[
e=\partial_U,\qquad
\partial_{V_i}\circ\partial_{V_j}=g_{ij}\partial_W,
\qquad
\partial_{V_i}\circ\partial_W
=\partial_W\circ\partial_W=0,
\]
where $g=(g_{ij})$ is a constant nondegenerate matrix. Suppose also that a continuous derivation $\mathcal G$ satisfies
\[
\mathcal G|_R=\delta,\qquad
\mathsf{Lie}_{\mathcal G}(\circ)=\circ,\qquad
\mathcal G\equiv\delta+U\partial_U-W\partial_W\pmod{F^1}.
\]
Then a formal coordinate change over $R$, equal to the identity modulo $F^1$, puts the product in the form
\begin{equation}\label{eq:negative-normal-product}
\begin{gathered}
e=\partial_U,\qquad
\partial_{V_i}\circ\partial_{V_j}=G_{ij}\partial_W,\\
\partial_{V_i}\circ\partial_W=0,
\qquad \partial_W\circ\partial_W=0.
\end{gathered}
\end{equation}
The matrix $G=(G_{ij})$ is invertible and satisfies $G\equiv g\pmod{F^1}$. Its entries are independent of $U, W$ and involve only curve classes $\beta$ with $\lambda_\beta=0$.

After extending coefficients to $\Omega:=\overline{\operatorname{Frac}(\mathcal O)}$, the tangent algebra is local of dimension $s+2$. Its radical has dimension $s+1$, its radical square has dimension one, and its radical cube is zero.
\end{lemma}

\begin{proof}
We first choose coordinates in which the unit is constant. We then simplify $\mathcal G$ sufficiently to read the possible product coefficients from the identity $\mathsf{Lie}_{\mathcal G}(\circ)=\circ$.

Since $e\equiv\partial_U\pmod{F^1}$, formal integration along $e$ gives a coordinate change, equal to the identity modulo $F^1$, for which $e=\partial_U$. The Hertling--Manin identity with both fields equal to $e$ gives $\mathsf{Lie}_e(\circ)=2\mathsf{Lie}_e(\circ)$, hence
\[
\mathsf{Lie}_{\partial_U}(\circ)=0.
\]
Thus every product coefficient is independent of $U$. Applying $\mathsf{Lie}_{\mathcal G}(\circ)=\circ$ to the unit also gives $[\mathcal G,e]=-e$. Consequently,
\[
\mathcal G=\mathcal G_0+N,
\qquad
\mathcal G_0:=\delta+U\partial_U-W\partial_W,
\]
where $N$ is $R$-linear, has coefficients in $F^1$, and is independent of $U$.

\medskip

To simplify $N$, we compute the effect of a coordinate change on each of its monomials. Assign the auxiliary weights
\[
\omega_U=1,\qquad \omega_{V_i}=0,
\qquad \omega_W=-1.
\]
These are weights for the present calculation, not the cohomological degrees used earlier. With $V^\nu:=\prod_iV_i^{\nu_i}$ and $y_j\in\{U,V_1,\ldots,V_s,W\}$, we have
\begin{equation}\label{eq:homological-eigenvalue}
[\mathcal G_0,q_S^\beta V^\nu W^k\partial_{y_j}]
=(\lambda_\beta-k-\omega_j)
q_S^\beta V^\nu W^k\partial_{y_j}.
\end{equation}
The three contributions are the action of $\delta$ on $q_S^\beta$, the action of $-W\partial_W$ on $W^k$, and the bracket with $\partial_{y_j}$.

Write $\mu:=\lambda_\beta-k-\omega_j$. If $\mu\neq0$, a term $a q_S^\beta V^\nu W^k\partial_{y_j}$ can be removed to its current $F$-order by the coordinate change generated by
\[
H:=\frac{a}{\mu}q_S^\beta V^\nu W^k\partial_{y_j}.
\]
Indeed, $[H,\mathcal G_0]=-a q_S^\beta V^\nu W^k\partial_{y_j}$. We call a term \emph{resonant} if $\mu=0$; such a term commutes with $\mathcal G_0$ and is retained.

At each order $r\geq1$, sum the indicated generators over all nonresonant monomials in $N$ modulo $F^{r+1}$, obtaining an $R$-linear derivation $H_r\in F^r$. This is well-defined coefficientwise as bounded projected degree leaves finitely many classes, and the remaining sums are formal power series. Conjugating derivations by $\exp(H_r)$ replaces $\mathcal G$ by $\exp(\operatorname{ad}H_r)\mathcal G$ and removes these terms modulo $F^{r+1}$. Brackets with $N$ lie in $F^{r+1}$, so lower orders remain unchanged.

Each change fixes $R$ and preserves $\partial_U$. Since $H_r$ raises the filtration, its exponential is well-defined, and the successive changes and their inverses converge. Renaming the resulting coordinates, we obtain
\begin{equation}\label{eq:resonant-normal-field}
\mathcal G=\mathcal G_0+N,
\qquad N\in F^1,
\qquad [\mathcal G_0,N]=0.
\end{equation}
The correction $N$ need not vanish.

\medskip

Let $M$ denote the multiplication tensor, namely the bilinear map $(A, B)\mapsto A\circ B$ in the new coordinates. The rescaling identity becomes
\begin{equation}\label{eq:product-weight-equation}
(\mathsf{Lie}_{\mathcal G_0}-1)M
=-\mathsf{Lie}_N M.
\end{equation}
Consider a term $q_S^\beta V^\nu W^k$ in the coefficient of $\partial_{y_a}$ in $\partial_{y_i}\circ\partial_{y_j}$.
Its contribution to $\mathsf{Lie}_{\mathcal G_0}M$ is multiplied by
\[
\lambda_\beta-k+\omega_i+\omega_j-\omega_a.
\]
We call this number its tensor weight. The leading product, obtained modulo $F^1$, has tensor weight one.

Because $[\mathcal G_0,N]=0$, the operator $\mathsf{Lie}_N$ preserves tensor weights. It also raises $F$-order, since $N$ has coefficients in $F^1$. Group the monomial terms of $M$ by their tensor weight and write $M_\mu$ for the terms of weight $\mu$. For $\mu\neq1$, equation~\eqref{eq:product-weight-equation} gives
\[
(\mu-1)M_\mu=-\mathsf{Lie}_N M_\mu.
\]
Initially $M_\mu\in F^1$. If $M_\mu\in F^r$, the right-hand side belongs to $F^{r+1}$, so the same holds for $M_\mu$, because $\mu-1$ is invertible. Induction and separatedness of the filtration imply $M_\mu=0$. Therefore \begin{equation}\label{eq:product-weight-one}
\mathsf{Lie}_{\mathcal G_0}M=M.
\end{equation}

We can now determine the possible products. Their coefficients are independent of $U$, so every coefficient monomial has weight $\lambda_\beta-k\leq0$. On the other hand, \eqref{eq:product-weight-one} requires the coefficient of $\partial_{y_a}$ in $\partial_{y_i}\circ\partial_{y_j}$ to have weight
\[
1-\omega_i-\omega_j+\omega_a.
\]
For products not involving the unit, these required weights are
\[
\begin{array}{c|ccc}
 &\partial_U&\partial_{V_a}&\partial_W\\ \hline
\partial_{V_i}\circ\partial_{V_j}&2&1&0\\
\partial_{V_i}\circ\partial_W&3&2&1\\
\partial_W\circ\partial_W&4&3&2
\end{array}
\]
Every positive entry forces the corresponding coefficient to vanish. Only the $\partial_W$ component of $\partial_{V_i}\circ\partial_{V_j}$ can remain. Its weight is zero, so $\lambda_\beta-k=0$. Since $\lambda_\beta\leq0$ and $k\geq0$, this requires $\lambda_\beta=0$ and $k=0$. This proves the form \eqref{eq:negative-normal-product} and the asserted dependence of $G_{ij}$. Their reduction is $g_{ij}$, so $G$ is invertible.

\medskip

Finally, let
\[
\mathfrak j:=
\bigoplus_{i=1}^s\mathcal O\partial_{V_i}
\oplus\mathcal O\partial_W.
\]
The displayed multiplication shows that $\mathfrak j$ is an ideal and that
\[
\mathfrak j^2=\mathcal O\partial_W,
\qquad \mathfrak j^3=0,
\qquad
\bigl(\operatorname{Der}^{\mathrm{cont}}_R(\mathcal O),\circ\bigr)
/\mathfrak j\simeq\mathcal O.
\]
For the equality in the square, choose $i,j$ with $g_{ij}\neq0$; the corresponding $G_{ij}$ is a unit. After extension to $\Omega$, the ideal $\mathfrak j\otimes_{\mathcal O}\Omega$ is nilpotent, and the quotient is the field $\Omega$. It is therefore the unique maximal ideal. Its dimension is $s+1$, its square has dimension one, and its cube is zero, as claimed.
\end{proof}

\subsection{Proof of Theorem~\ref{thm:conic-bundle-local-factor}}
 We may assume that $\Delta\ne\emptyset$ since otherwise $\Delta_S$ is not nef. Put $n:=\rho_S+2$. Since $\Delta_S$ is nef and the intersection matrix $(g_{ij})$ is nondegenerate, Lemma~\ref{lem:nef-nonsplitting} applies to the negative factor by \eqref{eq:negative-leading-field}. Extend coefficients to
\[
\Omega:=\overline{\operatorname{Frac}(\widehat{\mathcal O})},
\]
and denote the resulting two algebras by $B_+$ and $B_-$. The lemma gives an ideal
\[
\mathfrak m_-=
\bigoplus_{i=1}^{\rho_S}\Omega\partial_{V_i}
\oplus\Omega\partial_W
\subset B_-
\]
such that
\[
B_-/\mathfrak m_-\simeq\Omega,\qquad
\mathfrak m_-^2=\Omega\partial_W,\qquad
\mathfrak m_-^3=0.
\]
Here the equality for the square follows from the invertibility of the matrix $(G_{ij})$. Thus $B_-$ is local, with maximal ideal $\mathfrak m_-$, and
\[
\dim_\Omega B_-=n,\qquad
\dim_\Omega\mathfrak m_-=n-1,\qquad
\dim_\Omega\mathfrak m_-^2=1.
\]

It remains to check that its eigenvalue does not occur on $B_+$. Let $P_\pm(T)\in\widehat{\mathcal O}[T]$ be the characteristic polynomials of multiplication by the Euler vector field on the two lifted factors. The vertical formulas give, modulo $F^1$,
\[
\overline P_\pm(T)=(T-\lambda_\pm^{(0)})^n,
\qquad
\lambda_\pm^{(0)}=2q_f+u'\pm2q_fe^{z/2}.
\]
Consequently,
\[
\operatorname{Res}_T(\overline P_+,\overline P_-)
=\pm\bigl(4q_fe^{z/2}\bigr)^{n^2}
\]
is a unit. Completeness for the filtration implies that $\operatorname{Res}_T(P_+, P_-)$ is also a unit. Hence $P_+$ and $P_-$ have no common root over $\Omega$.

Since $B_-$ is local, multiplication by the Euler vector field has a single eigenvalue on $B_-$. Its algebraic multiplicity on $B_+\times B_-$ is therefore exactly $n$.

Finally, apply Section~\ref{sec:coefficient-comparison} to the injective homomorphism \eqref{eq:conic-coefficient-embedding}. The subsequent formal coordinate changes are invertible. Local-factor dimensions, radical filtrations, and algebraic multiplicities are unchanged after extension to algebraically closed fraction fields. Thus
\[
QH^{\mathrm{even}}(X)
\otimes_{\mathsf{Nov}^{\mathrm{red}}_X}\mathbb K_X
\]
has a local factor $\mathcal H_-$ with the asserted dimensions and radical filtration, and $\boldsymbol\kappa_{t'}$ has an eigenvalue of algebraic multiplicity $\rho_S+2$. Since $\rho_S\geq1$, this multiplicity is at least three. Proposition~\ref{thm:cai-bound} therefore implies that $X$ is irrational.

\section{Birational reduction and irrationality}
\label{sec:conic-bundle-birational-reduction}

Theorem~\ref{thm:conic-bundle-local-factor} requires $\Delta_S=2K_S+\Delta$ to be nef. We now show that, if
$|\Delta_S|\neq\emptyset$, there is a birational standard model whose corresponding class is nef. Here the hypothesis means that $\Delta_S$ is linearly equivalent to an effective divisor, possibly zero. This proves the irrationality implication in Shokurov's conjecture \cite[Conjecture~1.2]{Prokhorov}:

\begin{theorem}\label{thm:conic-bundle-effective}
Let $\pi: X\to S$ be a projective standard conic bundle over a smooth
projective rational surface, with discriminant $\Delta$. If
\[
|\Delta_S|\neq\emptyset,
\]
then $X$ is irrational.
\end{theorem}

First, recall from \cite[Sections~3.8--3.9]{Prokhorov} that for every decomposition $\Delta=A+B$ into effective divisors without common components, the intersection number $A\cdot B$
is even. Moreover, every connected component of $\Delta$ has arithmetic genus at least one. Consequently, every smooth rational component $\Gamma\subset\Delta$ satisfies
\begin{equation}\label{eq:discriminant-rational-component}
\Gamma\cdot(\Delta-\Gamma)\geq2.
\end{equation}
Indeed, this intersection number is even, and if it were zero, $\Gamma$ would itself be a connected component of $\Delta$.

The birational construction we need is \cite[Proposition~8.5]{Prokhorov}. If $\Gamma\subset S$ is a
$(-1)$-curve satisfying one of
\begin{equation}\label{eq:conic-contraction-cases}
\begin{aligned}
&\Gamma\not\subset\Delta,\qquad \Gamma\cdot\Delta=0,\\
&\Gamma\not\subset\Delta,\qquad \Gamma\cdot\Delta=1,\\
&\Gamma\subset\Delta,\qquad \Gamma\cdot(\Delta-\Gamma)=2,
\end{aligned}
\end{equation}
then its contraction $f:S\to S'$ lifts to a birational map $\varphi:X\dashrightarrow X'$ and a projective standard conic bundle $\pi':X'\to S'$ such that
\[
\pi'\circ\varphi=f\circ\pi,\qquad \Delta'=f_*\Delta.
\]
The construction gives a smooth total space $X'$. Since $\pi'$ is projective and $S'$ is projective, $X'$ is projective as well.

\begin{proposition}\label{prop:conic-nef-reduction}
Under the hypotheses of Theorem~\ref{thm:conic-bundle-effective},
there is a birational morphism $f:S\to S_{\mathrm{nef}}$ and a
projective standard conic bundle
\[
\pi_{\mathrm{nef}}:X_{\mathrm{nef}}\longrightarrow S_{\mathrm{nef}}
\]
with smooth total space, birational to $X$ over $f$, such that $S_{\mathrm{nef}}$ is a smooth projective rational surface and
\[
\Delta_{\mathrm{nef}}\neq\emptyset,\qquad
2K_{S_{\mathrm{nef}}}+\Delta_{\mathrm{nef}}
\quad\text{is nef and linearly equivalent to an effective divisor}.
\]
Here, $\Delta_{\mathrm{nef}}$ denotes the discriminant of $\pi_{\mathrm{nef}}$.
\end{proposition}
\begin{proof}
If $\Delta_S$ is already nef, no contraction is needed. Otherwise, choose an irreducible curve $\Gamma$ with $\Delta_S\cdot\Gamma<0$. Let $E\in|\Delta_S|$ be an effective representative. Since distinct irreducible curves on a smooth surface have nonnegative intersection, $\Gamma$ must occur in $E$. Writing $E=m\Gamma+E_0$, with $m>0$ and $\Gamma\not\subset\operatorname{Supp}(E_0)$, we obtain
\[
0>E\cdot\Gamma=m\Gamma^2+E_0\cdot\Gamma,
\]
so $\Gamma^2<0$. Write 
\[
\Gamma^2=-a,\qquad g=p_a(\Gamma),\qquad a\geq1,
\]
where $p_a$ denotes the arithmetic genus; adjunction gives $K_S\cdot\Gamma=2g-2+a$.

If $\Gamma\not\subset\Delta$, then
\[
0>\Delta_S\cdot\Gamma
=4g-4+2a+\Delta\cdot\Gamma.
\]
Since $g\geq0$ and $\Delta\cdot\Gamma\geq0$, this forces $g=0$, $a=1$, and $\Delta\cdot\Gamma\in\{0,1\}$.

If $\Gamma\subset\Delta$, put $b:=\Gamma\cdot(\Delta-\Gamma)\geq0$. Then
\[
0>\Delta_S\cdot\Gamma=4g-4+a+b.
\]
Thus, $g=0$, and hence $\Gamma$ is smooth and rational. By \eqref{eq:discriminant-rational-component}, $b\geq2$. The inequality $a+b<4$ therefore gives $a=1$ and $b=2$. In both cases, $\Gamma$ is a $(-1)$-curve satisfying one of
\eqref{eq:conic-contraction-cases}.

Let $f:S\to S'$ be the contraction of $\Gamma$, and apply \cite[Proposition~8.5]{Prokhorov}. Put $\Delta_{S'}:=2K_{S'}+\Delta'$. Since $K_S\sim f^*K_{S'}+\Gamma$, the three cases give
\begin{equation}\label{eq:adjoint-under-contraction}
\Delta_S\sim f^*\Delta_{S'}+
\begin{cases}
2\Gamma,&\Gamma\not\subset\Delta,\ \Gamma\cdot\Delta=0,\\
\Gamma,&\text{in the other two cases}.
\end{cases}
\end{equation}
To see this, $f^*\Delta'=\Delta$ in the first case and $f^*\Delta'=\Delta+\Gamma$ in the second. In the third, $\Delta-\Gamma$ meets $\Gamma$ transversely at two distinct points, because $\Delta$ is nodal. Its image has a node at $f(\Gamma)$, so
\[
f^*\Delta'=(\Delta-\Gamma)+2\Gamma=\Delta+\Gamma.
\]
Pushing forward an effective representative of $\Delta_S$ now shows that $|\Delta_{S'}|\neq\emptyset$.

Repeat this construction whenever the resulting adjoint class is not nef. Effectivity is preserved at each step, and every contraction lowers the Picard number of the base by one. The process therefore terminates at a smooth projective rational surface $S_{\mathrm{nef}}$ with nef adjoint class. The successive standard models remain smooth, projective, and birational to $X$.  Finally, the terminal discriminant cannot be empty. Otherwise, $|2K_{S_{\mathrm{nef}}}|\neq\emptyset$, contradicting
\[
H^0\bigl(S_{\mathrm{nef}},
\mathcal O_{S_{\mathrm{nef}}}(2K_{S_{\mathrm{nef}}})\bigr)=0,
\]
since every smooth projective rational surface has vanishing positive plurigenera.
\end{proof}

\begin{proof}[Proof of Theorem~\ref{thm:conic-bundle-effective}]
Proposition~\ref{prop:conic-nef-reduction} gives a projective standard conic bundle
$\pi_{\mathrm{nef}}:X_{\mathrm{nef}}\to S_{\mathrm{nef}}$ with nonempty discriminant and nef class
$2K_{S_{\mathrm{nef}}}+\Delta_{\mathrm{nef}}$. Theorem~\ref{thm:conic-bundle-local-factor} implies that $X_{\mathrm{nef}}$ is irrational. Since $X$ is birational to $X_{\mathrm{nef}}$, it is irrational as well.
\end{proof}

\medskip

Let $J(X)$ denote the intermediate Jacobian of $X$ equipped with its canonical principal polarization. As a complex torus, it is
\[
J(X):=\frac{H^3(X,\mathbb C)}
{F^2H^3(X,\mathbb C)+\operatorname{im}H^3(X,\mathbb Z)},
\qquad
F^2H^3(X,\mathbb C):=H^{3,0}(X)\oplus H^{2,1}(X).
\]
For the conic bundles over rational surfaces considered here, this is a principally polarized abelian variety. All isomorphisms below preserve the principal polarizations.

\begin{corollary}\label{cor:conic-rationality-criterion}
Let $\pi\colon X\to S$ be a standard conic bundle over a smooth projective rational surface, with discriminant $\Delta$.
Then $X$ is rational if and only if $|2K_S+\Delta|=\emptyset$ and the intermediate Jacobian $J(X)$ is isomorphic, as a principally polarized abelian variety, to a product of Jacobians of smooth projective curves.
\end{corollary}
\begin{proof}
If $X$ is rational, Theorem~\ref{thm:conic-bundle-effective} gives $|2K_S+\Delta|=\emptyset$, and the intermediate-Jacobian obstruction implies that $J(X)$ is a product of curve Jacobians as a principally polarized abelian variety \cite[Corollary~7.3]{Prokhorov}.

Conversely, assume both conditions. If $\Delta=\emptyset$, the vanishing of $\operatorname{Br}(S)$ gives a rational section, so $X$ is birational to $S\times\mathbb P^1$ and is rational \cite[Lemma~3.7]{Prokhorov}.

Suppose now that $\Delta\neq\emptyset$. By \cite[Proposition~8.6]{Prokhorov}, $X\to S$ is fiberwise birational to a standard conic bundle $X'\to S'$ with discriminant $\Delta'$, where either $S'$ is a Hirzebruch surface with $\Delta'\cdot F\leq3$ for a ruling fiber $F$, or $S'\simeq\mathbb P^2$. In the ruled case, $X'$ is rational by \cite[Proposition~5.6]{Prokhorov}. In the plane case, $J(X')$ is also a product of curve Jacobians as a principally polarized abelian variety. Indeed, this property is birationally invariant by the blowup formula and uniqueness of decomposition into indecomposable principally polarized factors \cite[Proposition~7.2]{Prokhorov}. Thus \cite[Theorem~9.1]{Prokhorov} implies that $X'$ is rational. Since $X$ is birational to $X'$, it is rational as well.
\end{proof}

\section{Nodal quartic double solids}
\label{sec:nodal-quartic-double-solids}

We now apply the conic-bundle calculation to double covers of $\mathbb P^3$.
By a \emph{nodal quartic double solid} we mean a double cover
$\tau:Y\to\mathbb P^3$ branched over a quartic surface
$B\subset\mathbb P^3$ whose only singularities are ordinary double
points; the smooth case is allowed. The nodes of $Y$ correspond
bijectively to those of $B$, and we write $N$ for their number.
Recall that $Y$ is \emph{$\mathbb Q$-factorial} if every Weil divisor
on $Y$ has a positive multiple that is Cartier.

We first construct a smooth standard conic-bundle model in the
$\mathbb Q$-factorial case and apply our quantum calculation.
Combining the resulting irrationality statement with the known
non-$\mathbb Q$-factorial cases then gives a rationality criterion
in terms of the nodes of $B$.

\subsection{The standard model and its quantum cohomology}

The geometric input for the quantum calculation is the following standard-model construction.
\begin{proposition}[{\cite[Theorem~4.2]{CPS2019}}]
\label{prop:quartic-standard-model}
Let $Y$ be a $\mathbb Q$-factorial nodal quartic double solid with
$N\geq1$ nodes. There is a projective standard conic bundle
$\pi:X\to S$, with smooth total space $X$ birational to $Y$, such
that $S\to\mathbb P^2$ is a composition of $N-1$ blowups of smooth
points, possibly infinitely near, and its discriminant satisfies
\[
\Delta\sim-2K_S.
\]
In particular,
\[
\rho_S=N,\qquad K_S^2=10-N,\qquad \Delta_S\sim0.
\]
\end{proposition}

\begin{corollary}\label{cor:quartic-quantum-irrationality}
Every $\mathbb Q$-factorial nodal quartic double solid is
irrational. For $N\geq1$, the standard model $X$ of
Proposition~\ref{prop:quartic-standard-model} has a local Artinian
factor $\mathcal H_-$ in
\[
QH^{\mathrm{even}}(X)
\otimes_{\mathsf{Nov}^{\mathrm{red}}_X}\mathbb K_X
\]
satisfying
\[
\begin{aligned}
\dim_{\mathbb K_X}\mathcal H_-&=N+2,&
\dim_{\mathbb K_X}\operatorname{rad}(\mathcal H_-)&=N+1,\\
\dim_{\mathbb K_X}\bigl(\operatorname{rad}(\mathcal H_-)^2\bigr)&=1,&
\operatorname{rad}(\mathcal H_-)^3&=0.
\end{aligned}
\]
The operator $\boldsymbol\kappa_{t'}$ on the full even quantum cohomology of $X$ has an eigenvalue of algebraic multiplicity exactly $N+2$.
\end{corollary}
\begin{proof}
For $N\geq1$, Proposition~\ref{prop:quartic-standard-model} gives $\Delta_S\sim0$ and $\rho_S=N$. Moreover, $\Delta$ is nonempty; otherwise $2K_S\sim0$, contradicting the rationality of $S$.
Theorem~\ref{thm:conic-bundle-local-factor} therefore gives the asserted local factor and multiplicity. Since $N+2\geq3$, Proposition~\ref{thm:cai-bound} excludes rationality of $X$, and hence of $Y$. For $N=0$, irrationality is classical \cite[Theorem~1.2 and the preceding discussion]{CPS2019}.
\end{proof}

\subsection{Rationality and the configuration of nodes}

The preceding corollary settles the $\mathbb Q$-factorial case.
Together with the results of Cheltsov--Przyjalkowski--Shramov,
it gives the following classification.

\begin{theorem}\label{thm:quartic-rationality-classification}
Let $Y$ be a nodal quartic double solid with $N$ nodes. Then
\[
Y\text{ is rational}
\quad\Longleftrightarrow\quad
Y\text{ is not }\mathbb Q\text{-factorial and }N\neq6.
\]
\end{theorem}
\begin{proof}
If $Y$ is rational, Corollary~\ref{cor:quartic-quantum-irrationality} shows that it is not $\mathbb Q$-factorial, while
\cite[Theorem~1.2]{CPS2019} excludes $N=6$. Conversely, if $Y$ is not $\mathbb Q$-factorial and $N\neq6$, rationality follows from \cite[Theorem~1.5]{CPS2019}.
\end{proof}

The criterion can also be read directly from the nodes of the branch surface. Put $\Sigma:=\operatorname{Sing}(B)$. By Clemens'
criterion \cite[Theorem~1.6]{CPS2019}, $Y$ is $\mathbb Q$-factorial if and only if the evaluation map
\[
\operatorname{ev}_\Sigma:
H^0(\mathbb P^3,\mathcal O_{\mathbb P^3}(2))
\longrightarrow H^0(\Sigma,\mathcal O_\Sigma(2))
\]
is surjective. Write the nodes as
$P_j=[z_{j0}:z_{j1}:z_{j2}:z_{j3}]$ and choose nonzero vector
representatives $(z_{j0},\ldots,z_{j3})$. Using the ten quadratic
monomials $z_a z_b$, $0\leq a\leq b\leq3$, as a basis, the map
$\operatorname{ev}_\Sigma$ is represented by the $N\times10$ matrix
\[
(M_\Sigma)_{j,(a,b)}:=z_{ja}z_{jb}.
\]
Changing a representative rescales the corresponding row by a nonzero scalar, so the rank is independent of these choices.
Theorem~\ref{thm:quartic-rationality-classification} therefore gives
\[
Y\text{ is rational}
\quad\Longleftrightarrow\quad
\operatorname{rank}_{\mathbb C}M_\Sigma<N
\ \text{ and }\ N\neq6.
\]
In particular, every nodal quartic double solid with $N\geq11$ is rational, since $M_\Sigma$ has only ten columns.

\end{document}